\documentclass[12pt,amstex]{amsart}
\usepackage{amssymb,amsthm,amsmath,fullpage}
\usepackage{hyperref}

\DeclareFontFamily{OMS}{rsfs}{\skewchar\font'60}
\DeclareFontShape{OMS}{rsfs}{m}{n}{<-5>rsfs5 <5-7>rsfs7 <7->rsfs10 }{}
\DeclareSymbolFont{rsfs}{OMS}{rsfs}{m}{n}
\DeclareSymbolFontAlphabet{\scr}{rsfs}

\DeclareFontFamily{OT1}{pzc}{}
\DeclareFontShape{OT1}{pzc}{m}{it}%
             {<-> s * [1.100] pzcmi7t}{}
\DeclareMathAlphabet{\mathchanc}{OT1}{pzc}%
                                 {m}{it}

\def\coh#1.#2.#3.{H^{#1}(#2,#3)}
\def\ring#1.{\scr O_{#1}}

\makeatletter
\renewcommand\subsection{
  \renewcommand{\sfdefault}{pag}
  \@startsection{subsection}%
  {2}{0pt}{-\baselineskip}{.2\baselineskip}{\raggedright
    \sffamily\itshape\small
  }}
\renewcommand\section{
  \renewcommand{\sfdefault}{phv}
  \@startsection{section} %
  {1}{0pt}{\baselineskip}{.2\baselineskip}{\centering
    \sffamily
    \scshape
}}
\makeatother

\newcommand{\mydot}{{{\,\begin{picture}(1,1)(-1,-2)\circle*{2}\end{picture}\ }}}
\newcommand{\DuBois}[1]{{\underline \Omega {}^0_{#1}}}
\newcommand{\FullDuBois}[1]{{\underline \Omega {}^{\mydot}_{#1}}}
\newcommand{\qis}{\simeq_{\text{qis}}}
\newcommand{\blank}{\underline{\hskip 10pt}}
\newcommand{\tensor}{\otimes}
\newcommand{\into}{\hookrightarrow}
\newcommand{\sI}{\scr{I}}
\newcommand{\sL}{\scr{L}}
\newcommand{\sK}{\scr{K}}
\newcommand{\sO}{\scr{O}}
\newcommand{\sF}{\scr{F}}
\newcommand{\sG}{\scr{G}}
\newcommand{\sA}{\scr{A}}

\newcommand{\OO}{\sO} 
\newcommand{\tld}{\widetilde }
\newcommand{\bC}{\mathbb{C}}

\newcommand{\bZ}{\mathbb{Z}}
\newcommand{\bH}{\mathbb{H}}
\newcommand{\bQ}{\mathbb{Q}}
\newcommand{\bA}{\mathbb{A}}
\newcommand{\myR}{\mcR\!} 
\newcommand{\mco}{\mathchanc{o}}
\newcommand{\mcm}{\mathchanc{m}}
\newcommand{\mcH}{\mathchanc{H}}

\newcommand{\mcR}{\mathchanc{R}}
\newcommand{\omegaplus}[1]{\underline\Omega_{{#1}}^{\times}}

\DeclareMathOperator{\coherent}{{coh}}
\DeclareMathOperator{\quasicoherent}{{qcoh}}
\DeclareMathOperator{\sHom}{\!\mcH\mco\mcm}

\DeclareMathOperator{\exc}{{Ex}}
\DeclareMathOperator{\Supp}{{Supp}}
\DeclareMathOperator{\Sing}{{Sing}}
\DeclareMathOperator{\divisor}{{div}}

\DeclareMathOperator{\Ann}{{Ann}}
\DeclareMathOperator{\red}{red}
\DeclareMathOperator{\codim}{codim}
\DeclareMathOperator{\supp}{{supp}}

\input xy
\xyoption{all}

\newtheoremstyle{dubois}{3pt}{3pt}%
     {\itshape}
     {}
     {\bfseries}
     {.}
     {.5em}
     {\thmname{#1}\thmnumber{ #2}\thmnote{ \rm #3}}
\newtheoremstyle{dubois-sf}{3pt}{3pt}%
     {\itshape}
     {}
     {\sffamily}
     {.}
     {.5em}
     {\thmname{#1}\thmnumber{ #2}\thmnote{ \rm #3}}
\newtheoremstyle{dubois-sc}{3pt}{3pt}%
     {\slshape}
     {}
     {\scshape}
     {.}
     {.5em}
     {\thmname{#1}\thmnumber{ #2}\thmnote{ \rm #3}}
\newtheoremstyle{dubois-remark}{3pt}{3pt}%
     {} 
     {}
     {\scshape}
     {.}
     {.5em}
     {\thmname{#1}\thmnumber{ #2}\thmnote{ \rm #3}}
\newtheoremstyle{dubois-def}{3pt}{3pt}%
     {}
     {}
     {\bfseries}
     {.}
     {.5em}
     {\thmname{#1}\thmnumber{ #2}\thmnote{ \rm #3}}
\newtheoremstyle{dubois-claim}{3pt}{3pt}%
     {}
     {}
     {\slshape}
     {:}
     {.5em}
     {\thmname{#1}\thmnote{ \rm #3}}

\theoremstyle{dubois}
\newtheorem{theorem}{Theorem}[section]
\newtheorem{lemma}[theorem]{Lemma}
\newtheorem{proposition}[theorem]{Proposition}
\newtheorem{corollary}[theorem]{Corollary}

\theoremstyle{dubois-def}
\newtheorem{definition}[theorem]{Definition}

\newtheorem{claim}[equation]{\sl Claim}

\theoremstyle{dubois-claim}
\newtheorem{claim*}{Claim}

\theoremstyle{dubois-remark}
\newtheorem{remark}[theorem]{Remark}

\theoremstyle{dubois-sc}
\newtheorem{setup}[theorem]{Setup} 
\newtheorem{notation}[theorem]{Notation}
\newtheorem{convention}[theorem]{Convention}

\numberwithin{equation}{theorem}

\newtheorem{proclaim-special}[theorem]{\specialthmname}

\begin{document}

\title{
  The canonical sheaf of Du~Bois singularities}

\author{S\'andor J Kov\'acs, Karl Schwede and Karen E.\ Smith}

\address{\vskip-.8cm\ \newline \noindent S\'andor J.\ Kov\'acs: \sf University of
  Washington, Department of Mathematics, Seattle, WA 98195, USA}
\email{kovacs@math.washington.edu}

\address{\vskip -.65cm \noindent Karl E.\ Schwede: \sf Department of Mathematics,
  University of Michigan.  Ann Arbor, Michigan 48109-1109}
\email{kschwede@umich.edu}

\address{\vskip -.65cm \noindent Karen E.\ Smith: \sf Department of Mathematics,
  University of Michigan.  Ann Arbor, Michigan 48109-1109}
\email{kesmith@umich.edu}

\date{\today}

\thanks{The first named author was supported in part by NSF Grants DMS-0554697 and
  DMS-0856185, and the Craig McKibben and Sarah Merner Endowed Professorship in
  Mathematics.}

\thanks{The second named author was partially supported by RTG grant number 0502170
  and by a National Science Foundation Postdoctoral Research Fellowship.}

\thanks{The third named author was partially supported by NSF Grant
  DMS-0500823}

\maketitle

\centerline{\sf In memoriam Juha Heinonen}

\begin{abstract}
  We prove that a Cohen-Macaulay normal variety $X$ has Du~Bois singularities if and
  only if $\pi_*\omega_{X'}(G) \simeq \omega_X$ for a log resolution $\pi: X'
  \rightarrow X$, where $G$ is the reduced exceptional divisor of $\pi$. Many basic
  theorems about Du~Bois singularities become transparent using this characterization
  (including the fact that Cohen-Macaulay log canonical singularities are Du~Bois).
  We also give a straightforward and self-contained proof that (generalizations of)
  semi-log-canonical singularities are Du~Bois, in the Cohen-Macaulay case.  It also
  follows that the Kodaira vanishing theorem holds for semi-log-canonical varieties
  and that Cohen-Macaulay semi-log-canonical singularities are cohomologically
  insignificant in the sense of Dolgachev.
\end{abstract}

\section{Introduction}

Consider a complex algebraic variety $X$.  If $X$ is smooth and projective, its De
Rham complex plays a fundamental role in understanding the geometry of $X$. When $X$
is singular, an analog of the De Rham complex, the \emph{Deligne-Du~Bois complex}
plays a similar role.  Based on Deligne's theory of mixed Hodge structures, Du~Bois
defined a filtered complex of $\sO_X$-modules, denoted by $\FullDuBois{X}$, that
agrees with the algebraic De Rham complex in a neighborhood of each smooth point and,
like the De Rham complex on smooth varieties, its {analytization} provides a
resolution of the sheaf of locally constant functions on $X$ \cite{DuBoisMain}.

Du~Bois observed that an important class of singularities are those for which
$\DuBois{X}$, the zeroth graded piece of the filtered complex $\FullDuBois{X}$, takes
a particularly simple form (see Discussion \ref{SubsectionDB}). He pointed out that
singularities satisfying this condition enjoy some of the nice Hodge-theoretic
properties of smooth varieties.  Dubbed {\it Du~Bois singularities \/} by Steenbrink,
these singularities have been promoted by Koll\'ar as a natural setting for vanishing
theorems; see, for example, \cite[Ch 12]{KollarShafarevich}.  Since the 1980's,
Steenbrink, Koll\'ar, Ishii, Saito and many others have investigated the relationship
between Du~Bois (or \emph{DB}) singularities and better known singularities in
algebraic geometry, such as rational singularities and log canonical singularities.
Because of the difficulties in defining and understanding the Deligne-Du~Bois complex
$\DuBois X$, many basic features of DB singularities have been slow to reveal
themselves or have remained obscure.  The purpose of this paper is to prove a simple
characterization of DB singularities in the Cohen-Macaulay case, making many of
their properties and their relationship to other singularities transparent.

Let $\pi: \tld X \rightarrow X $ be a log resolution of a normal complex variety $X$,
and denote by $G$ the reduced exceptional divisor of $\pi$. By
Lemma~\ref{PushdownInjects} there exists a natural inclusion $\pi_*\omega_{\tld X}(G)
\hookrightarrow \omega_X$.
%
%
The main foundational result of this article is the following multiplier ideal-like
criterion for DB singularities:

\begin{theorem}[{(=Theorem~\ref{Kempf})}]\label{thm:main}
  Suppose that $X$ is normal and Cohen-Macaulay.  Let $\pi : X' \rightarrow X$ be any
  log resolution, and denote the reduced exceptional divisor of $\pi$ by $G$.  Then
  $X$ has DB singularities if and only if $\pi_* \omega_{X'}(G) \simeq
  \omega_X$.
\end{theorem}

Theorem~\ref{thm:main} is analogous to the following well-known criterion for
rational singularities due to Kempf: if $X$ is normal and Cohen-Macaulay, then $X$
has rational singularities if and only if the natural inclusion $\pi_*\omega_{\tld X}
\hookrightarrow \omega_X$ is an isomorphism. In particular, Theorem~\ref{thm:main}
immediately implies that rational singularities are Du~Bois, a statement that had
been conjectured by Steenbrink in \cite{SteenbrinkMixed} and later proved by Koll\'ar
\cite{KollarShafarevich} in the projective case and finally by Kov\'acs
\cite{KovacsDuBoisLC1}, and also independently by Saito \cite{SaitoMixedHodge} in
general.   Another immediate corollary is that normal quasi-Gorenstein DB
singularities are log canonical; see Section~\ref{SectionProofOfTheMainTheorem} for a
complete discussion.

In addition, this criterion shows that CM DB singularities relate to rational
singularities very much like log canonical singularities relate to (kawamata) log
terminal singularities. This has been a general belief all along, but we feel that
the criterion in Theorem~\ref{thm:main} supports this belief more than anything else
previously known.

A long-standing conjecture of Koll\'ar's predicts that log canonical singularities
are Du~Bois.  Using Theorem~\ref{thm:main}, it is easy to see that Koll\'ar's
conjecture holds in the Cohen-Macaulay case:

\begin{theorem}[{(=Theorem~\ref{CMandLCImpliesDuBois})}]
  \label{thm:CMandLCImpliesDuBois}
  Suppose that $X$ is normal and Cohen Macaulay, and that $\Delta $ is
  an effective $\bQ$-divisor on $X$ such that $K_X + \Delta$ is
  $\bQ$-Cartier.  If $(X, \Delta)$ is log canonical, then $X$ has Du
  Bois singularities.
\end{theorem}

In fact, we prove the stronger result that Cohen-Macaulay semi-log canonical
singularities are Du~Bois, see Theorem~\ref{CMslcImpliesDuBois}, (even more, we prove
that a generalization of semi-log canonical singularities are Du~Bois).  This is
based on a technical generalization of certain aspects of Theorem~\ref{thm:main} to
the non-normal case, treated in Section~\ref{SectionNonnormalCase}.  Many special
cases of Koll\'ar's conjecture had been known, including the isolated singularity
case \cite{IshiiIsolatedGorenstein,IshiiDuBoisSurfaceSings,IshiiIsolatedQGorenstein},
the Cohen-Macaulay case when the singular set is not too big see also
\cite{KovacsDuBoisLC1}, and the local complete intersection case
\cite{SchwedeEasyCharacterization}.

Very recently, Koll\'ar's conjecture that log canonical singularities are Du~Bois,
has been verified by Koll\'ar and the first named author, using recent advances in
the Minimal Model Program \cite{KK09}. In particular, there is now an independent and
more general result proving Theorem~\ref{thm:CMandLCImpliesDuBois}.  However, there
are several reasons why Theorem~\ref{thm:main} is still interesting (besides being
the first general result of this kind). Koll\'ar and Kov\'acs also prove that the
condition of being Cohen-Macaulay is constant in DB families. This means that a
stable smoothable variety is necessarily Cohen-Macaulay, and hence the above
condition is applicable.

Furthermore, Theorem~\ref{keySlcDuBoisTheorem} and Theorem~\ref{CMslcImpliesDuBois},
apply to non-normal singularities. These are the best results currently known. Notice
that one of the main applications of DB singularities is to moduli theory and
that is an arena where non-normal singularities may not be easily dismissed. Already
for degenerations of curves one must deal with non-normal singularities. In
particular, Koll\'ar's conjecture is actually important in the non-normal case, that
is, we want to know that semi-log canonical singularities are Du~Bois.

Another immediate corollary, again predicted by Koll\'ar, is that the Kodaira
Vanishing Theorem holds for generalizations of (semi-)log canonical varieties.
Fujino recently gave another proof of a closely related theorem using techniques of
Ambro; see \cite[Corollary 5.11]{FujinoVanishingAndInjectivity}.

\begin{corollary}[{(=Corollary~\ref{KodairaVanishingForLC})}]
  Kodaira vanishing holds for Cohen-Macaulay weakly semi-log canonical varieties.  In
  particular, let $(X,\Delta)$ be a projective Cohen-Macaulay weakly semi-log
  canonical pair and $\sL$ an ample line bundle on $X$. Then $\coh i.X.\sL^{-1}.=0$
  for $i<\dim X$.
\end{corollary}

Of course, if $X$ is not Cohen-Macaulay, Kodaira vanishing in the
above form necessarily fails.  But the Cohen-Macaulay condition is not
sufficient for Kodaira vanishing.  Examples show that some further
restriction on the singularities is needed; see \cite[Section
2]{ArapuraJaffeKodairaVanishingForSingular}.  In some sense, this
is the most general form of the classical Kodaira vanishing theorem
(that is, $H^i(X, \sL^{-1}) = 0$ for $\sL$ ample and $i < \dim X$)
that could be hoped for.

We are also able to obtain some nice Hodge-theoretic properties for semi-log
canonical singularities.  In particular, we are able to show that Cohen-Macaulay
semi-log canonical singularities are cohomologically insignificant in the sense of
Dolgachev \cite{DolgachevCohomologicallyInsignificant}; see
Theorem~\ref{CMslcImpliesInsignificant}.  This fact has useful applications in the
construction of compact moduli spaces of stable surfaces and higher dimensional
varieties.


\section{Preliminaries}

In this section we will define the notion of log canonical, as well as DB
singularities and state the forms of duality we will use.  Throughout this paper, a
\emph{scheme} will always be assumed to be separated and noetherian of essentially
finite type over $\bC$.  By a \emph{variety}, we mean a reduced separated noetherian
pure-dimensional scheme of finite type over $\bC$.  Note that a variety may have
several irreducible components.  All varieties and schemes will be assumed to be
quasi-projective. The purpose of this assumption is to guarantee that these varieties
are embedded in smooth schemes.  Note that in the end this hypothesis is harmless
because implications between various types of singularities are local questions, thus
the varieties may assumed to be quasi-projective.

We will use the following notation: For a functor $\Phi$, $\myR\Phi$ denotes its
derived functor on the (appropriate) derived category and $\myR^i\Phi:=
h^i\circ\myR\Phi$ where $h^i(C^\mydot)$ is the cohomology of the complex $C^\mydot$
at the $i^\text{th}$ term.  Similarly, $\bH^i_Z := h^i \circ \myR \Gamma_Z$ where
$\Gamma_Z$ is the functor of cohomology with supports along a subscheme $Z$.
Finally, $\sHom$ stands for the sheaf-Hom functor.

Let $\alpha:Y\to Z$ be a birational morphism and $\Delta \subseteq Z$ a $\bQ$-divisor.
Then $\alpha^{-1}_*\Delta$ will denote the proper transform of $\Delta$ on $Y$.




\subsection{Log Canonical Singularities} \label{SubsectionLCPreliminaries}

Let $X$ be a normal irreducible variety of pure dimension $d$. The
canonical sheaf $\omega_X$ of $X$ is the unique reflexive
$\OO_X$-module agreeing with the sheaf of regular differential
$d$-forms $\bigwedge^d\Omega_{X/\bC}$ on the smooth locus of $X$. A
canonical divisor is any member $K_X$ of the (Weil) divisor class
corresponding to $\omega_X$. See \ref{SubsectionDuality} for the
definition of the canonical sheaf on non-normal varieties.

A ($\bQ$-)Weil divisor $D$ is said to be $\bQ$-Cartier if, for some non-zero integer
$r$, the $\bZ$-divisor $rD$ is Cartier, meaning that it is given
locally as the divisor of some rational function on $X$. For such a
divisor, we can define the pullback $\pi^*D$, under any dominant
morphism $\pi$, to be the $\bQ$-divisor $\frac{1}{r}\pi^*(rD).$ We say
that $X$ is $\bQ$-Gorenstein if $K_X$ is $\bQ$-Cartier.

Now consider a pair $(X, \Delta),$ where $\Delta$ is an effective
$\bQ$-divisor such that $K_X + \Delta$ is $\bQ$-Cartier. In this case,
there exists a {\it log resolution\/} of the pair; that is, a proper
birational morphism from a smooth variety
$$
\pi: \tld X \rightarrow X
$$
such that $\exc(\pi)$ is a divisor and furthermore the set $\pi^{-1}(\Delta) \cup
\exc(\pi)$ is a divisor with simple normal crossing support.  Here $\exc(\pi)$
denotes the exceptional set of $\pi$.  Let $\tld\Delta $ denote the birational (or
proper, or strict) transform of $\Delta$ on $\tld X$, often denoted by
$\pi^{-1}_*\Delta$. Then there is a numerical equivalence of divisors
$$
K_{\tld X} + \tld \Delta - \pi^{*}(K_X + \Delta) \equiv \sum a_i E_i
$$
where the $E_i$ are the exceptional divisors of $\pi$ and the $a_i$
are some uniquely determined rational numbers. We can now define:

\begin{definition}\label{LCdef}
  The pair $(X, \Delta)$ is called \emph{log canonical} if
  $a_i \geq -1$ for all $i$.
\end{definition}

Definition~\ref{LCdef} is independent of the choice of log resolution.
For this and other details about log resolutions, log pairs, and log
canonical singularities see, for example, \cite{KollarMori}.

\subsection{Du~Bois Singularities}\label{SubsectionDB}

As mentioned in the introduction, DB singularities are defined using a fairly
complicated filtered complex $\FullDuBois{X}$,
which plays the role of the De Rham complex for singular varieties. It follows from
the construction that there is a natural map (in the derived category of
$\OO_X$-modules)
$$
\OO_X \rightarrow \DuBois{X},
$$
where $\DuBois{X}$ denotes the zeroth graded complex of the filtered complex
$\FullDuBois{X}$.  By definition, $X$ has {\it Du~Bois (or DB) singularities} if this
map is a quasi-isomorphism.  For a careful development of this point of view see
\cite{DuBoisMain}, \cite{GNPP}, or \cite{SteenbrinkgVanishing}.  However, a recent
result of the second author in \cite{SchwedeEasyCharacterization}
provides an alternate definition of DuBois singularities, which we now
review here (also compare with \cite{EsnaultHodgeTypeOfSubvarieties90}).

First, since the transverse union of two smooth varieties of different
dimensions is an important example of a DB singularity, we must
leave the world of irreducible (and even that of equidimensional)
varieties, and instead consider reduced schemes of finite type over
$\bC$.  Recall that a reduced subscheme $X$ of a smooth ambient
variety $Y$ is said to have \emph{normal crossings} if at every point
of $x$ there exists a regular system of parameters for $\OO_{Y, x}$
such that $X$ is defined by some monomials in these parameters.

Now suppose that $X$ is a reduced closed subscheme of a smooth ambient
$Y$.  Let $\pi : \tld Y \rightarrow Y$ be a proper birational morphism
such that:
\begin{itemize}
\item[(a)]  $\pi$ is an isomorphism outside $X$,
\item[(b)]  $\tld Y$ is smooth,
\item[(c)] $\pi^{-1}(X)$ has normal crossings (even though it may not
  be equidimensional)
\end{itemize} Such morphisms always exist by Hironaka's theorem; for
example, $\pi$ could be an embedded resolution of singularities for $X
\subset Y$, or a log resolution of the pair $(Y, X)$.  The conditions
(a)-(c) can be relaxed somewhat; see \cite[Proposition
2.20]{SchwedeTakagiRationalPairs} for a way to relax condition (a) and
\cite{SchwedeEasyCharacterization} for further discussion.  We set
$\overline X$ to be the reduced preimage of $X$ in $\tld Y$, and again
emphasize that $\overline X$ may not be equidimensional.  Under these
conditions, Schwede shows that the object $\myR \pi_* \OO_{\overline
  X}$ is naturally quasi-isomorphic to the object $\DuBois{X}$ defined
in \cite{DuBoisMain}. This leads to the following definition,
equivalent to Steenbrink's original definition:

\begin{definition}
  \label{DBdef}
  We say that $X$ has \emph{Du~Bois (or DB) singularities} if the natural map $\OO_X
  \rightarrow \myR \pi_* \OO_{\overline X}$ is a quasi-isomorphism.
\end{definition}

The fact that Definition~\ref{DBdef} is independent of the choice of
embedding and of the choice of resolution simply follows from the fact
that $\myR \pi_* \OO_{\overline X}$ is quasi-isomorphic to the
well-defined object $\DuBois{X}$.  See \cite{SchwedeEasyCharacterization} for details (and again, compare with \cite{EsnaultHodgeTypeOfSubvarieties90}).

\subsection{Dualizing Complexes and Duality}\label{SubsectionDuality}


For the convenience of the reader, we briefly review the form of
duality we will use.

Associated to every quasi-projective scheme $X$ of dimension $d$ and of finite type
over $\bC$ there exists a (normalized) dualizing complex $\omega_X^{\mydot} \in
D^b_{\coherent}(X)$, unique up to quasi-isomorphism.  To construct
$\omega_X^{\mydot}$ concretely, one may proceed as follows.  For a smooth irreducible
variety $Y$ of dimension $n$, take $\omega_Y^{\mydot}=\omega_Y[n]$, the complex that
has the canonical module $\det\Omega_{Y/\bC}$ in degree $-n$ and the zero module in
all the other spots.  Now, whenever $X \subseteq Y$ is a closed embedding of schemes
of finite type over a field, we have
\begin{equation*}
  \omega_X^{\mydot} =  \myR \sHom_{Y}^{\mydot}(\OO_X, \omega_Y^{\mydot}).
\end{equation*}
Because every quasi-projective scheme of finite type over $\bC$ embeds
in a smooth variety, this determines the dualizing complex on any such
finite type scheme.  Alternatively, one can also define
$\omega_X^{\mydot}$ as $h^! \bC$ where $h : X \rightarrow \bC$ is the structural
morphism.  See \cite{HartshorneResidues} and
\cite{ConradGDualityAndBaseChange} for details.

For an equidimensional $X$ admitting a dualizing complex one may
define the canonical sheaf to be the $\OO_X$-module $h^{-\dim
  X}(\omega_X^{\mydot})$, denoted by $\omega_X$. If $X$ is
Cohen-Macaulay, the dualizing complex turns out to be exact at all
other spots (just as in the case of smooth varieties above).  In
particular, for a Cohen-Macaulay variety of pure dimension $d$, we
have $\omega_X^{\mydot} = \omega_X[d]$.
%
If $X$ is normal and irreducible, $h^{-\dim X}(\omega_X^{\mydot})$
agrees with the canonical sheaf $\omega_X$ as defined in
\ref{SubsectionLCPreliminaries}, so there is no ambiguity of
terminology. More generally, if $X$ is Gorenstein in {co}dimension
%
%
one and satisfies Serre's $S_2$ condition, the canonical module $h^{-\dim
  X}(\omega_X^{\mydot}) = \omega_X$ is a rank one reflexive sheaf, and so corresponds
to a Weil divisor class.

A very important tool we need is Grothendieck duality:

\begin{theorem}\cite[III.11.1, VII.3.4]{HartshorneResidues}
  \label{GD}
  Let $f : X \rightarrow Y$ be a proper morphism between finite
  dimensional noetherian schemes.  Suppose that both $X$ and $Y$ admit
  dualizing complexes and that $\sF^{\mydot} \in
  D^{-}_{\quasicoherent}(X)$.  Then the duality morphism
  \begin{equation*}
    \myR f_* \myR \sHom_{X}^{\mydot}(\sF^{\mydot}, \omega_X^{\mydot})
    \to \myR \sHom_{Y}^{\mydot}(\myR f_* \sF^{\mydot}, \omega_Y^{\mydot}).
  \end{equation*}
  is a quasi-isomorphism.
\end{theorem}

\begin{remark}
In the previous theorem, $\omega_X^{\mydot}$ should be thought of as $f^! \omega_Y^{\mydot}$.
\end{remark}






\section{Proof of the Main Theorem}
\label{SectionProofOfTheMainTheorem}
In this section we prove our main result, a simple new
characterization of DB singularities in the normal Cohen-Macaulay case.

\begin{theorem}\label{Kempf}
  Suppose that $X$ is a normal Cohen-Macaulay variety. Let $\varrho : X' \rightarrow
  X$ be any log resolution and denote the reduced exceptional divisor of $\varrho$ by
  $G$.  Then $X$ has DB singularities if and only if $\varrho_* \omega_{X'}(G)
  \simeq \omega_X$.
\end{theorem}

The proof of this theorem will take most of the present section.  We will first show
that it is true for a special choice of log resolutions in
Corollary~\ref{CMNormalSurjectiveImpliesDuBois}. Then, in
Lemma~\ref{lem:independent}, we complete the proof by showing that the statement is
independent of the choice of the log resolution.

\smallskip
\noindent {\sf \small The following notation will be fixed for the rest of
  this section:}
\smallskip

  \begin{setup}
    \label{setupNormal}
    Let $X$ be a reduced equidimensional scheme of finite type over $\bC$ embedded in
    a smooth variety $Y$, and let $\Sigma$ denote a closed subscheme $\Sigma
    \subsetneqq X$ that contains the singular locus of $X$.  Assume that no
    irreducible component of $X$ is a hypersurface
    i.e., the codimension of every irreducible component of $X$ in $Y$ is at least
    two.
    Fix an embedded resolution $\pi : \tld Y \rightarrow Y$ of $X$ in $Y$, and let
    $\tld X$ denote the strict transform of $X$ on $\tld Y$. Further assume that
    \begin{itemize}
    \item[(i)] $\pi$ is an isomorphism over $X\setminus \Sigma$,
      and
    \item[(ii)] $\pi^{-1}(\Sigma)$ is a simple normal crossing divisor of $\tld Y$
      that intersects $\tld X$ in a simple normal crossing divisor of $\tld X$.
    \end{itemize}

    Let $E$ denote the reduced preimage of $\Sigma$ in $\tld Y$, and $\overline X$
    the reduced pre-image of $X$ in $\tld Y$. We emphasize that $\overline X$ is {\it
      not} equidimensional; in fact, $\overline X$ is the transverse union of the
    smooth variety $\tld X$ and the normal crossing divisor $E$.  We will frequently
    abuse notation and use $\pi$ to denote $\pi|_{\overline X}$.  Note that
    $\pi|_{\overline X}$ is a projective (respectively proper) morphism as long as
    $\pi$ is. Finally, recall that by \cite{SchwedeEasyCharacterization}
    $R\pi_*\sO_{\overline X}\qis\DuBois X$ and $R\pi_*\sO_{E}\qis\DuBois\Sigma$.
  \end{setup}

\smallskip

The outline of the proof of Theorem~\ref{Kempf} goes as follows.  Using the
Grothendieck dual form of Schwede's characterization of DB singularities stated
in Lemma~\ref{Karldual} below, it is clear that a reduced Cohen-Macaulay scheme $X$
of dimension $d$ has DB singularities if and only if
\begin{itemize}
\item[(i)] $\myR^i \pi_* \omega^{\mydot}_{\overline X} = 0$ \,\,\, for $i \neq -d$,
  and
\item [(ii)] $\myR^{-d} \pi_* \omega^{\mydot}_{\overline X} = \omega_X$.
\end{itemize}
Our main technical statement is Theorem~\ref{top}, which implies that for any normal
variety of dimension $d$, the sheaf $\myR^{-d} \pi_* \omega^{\mydot}_{\overline X}$
can be identified with $\pi_*\omega_{\tld X}(G)$.  This is proven by comparing the
dualizing complexes for $\overline X$, $\tld X$ and $E$ via the dual of the short
exact sequence
$$
0 \rightarrow \OO_{\tld X}(-E|_{\tld X}) \rightarrow \OO_{\overline X}
\rightarrow \OO_E \rightarrow 0.
$$
On the other hand, the vanishing statement of (i) follows from a reinterpretation of
a result of the first author \cite{KovacsDuBoisLC1} by the second author
\cite{SchwedeFInjectiveAreDuBois}.

We first state the following dual form of Schwede's characterization of DB
singularities.

\begin{lemma}
  \label{Karldual}
  $X$ has DB singularities if and only if the natural map
  $$
  \myR \pi_* \omega^{\mydot}_{\overline X} \rightarrow \omega_X^{\mydot}
  $$
  is a quasi-isomorphism.
\end{lemma}

\begin{proof}
  The result follows from Definition~\ref{DBdef} via a standard application of
  Grothendieck duality.
\end{proof}

Before beginning the proof of Theorem~\ref{Kempf}, we would like to make the
following suggestive observation.  If $X$ has pure dimension $d$ and $Y$ has dimension
$n$, then
$$
\xymatrix@R=0pt{%
  h^{-d}(\omega_{\overline X}^{\mydot}) \cong
  \omega_{\tld X}(E|_{\tld X}),\\
  h^{-n+1}(\omega_{\overline X}^{\mydot}) \cong \omega_E \text{, and} \\
  h^{i}(\omega_{\overline X}^{\mydot}) = 0 \text{ for $i$ not equal to
    $-n+1$ or $-d$.}  }
$$

To see this, note that $\tld X$ and $E$ have normal crossings and that
there exists a short exact sequence,
$$
0 \rightarrow \OO_{\tld X}(-E|_{\tld X}) \rightarrow
\OO_{\overline X} \rightarrow \OO_{E} \rightarrow 0.
$$
Next, dualize this sequence by applying $\myR \sHom_{\tld Y}^{\mydot}(\blank,
\omega_{\tld Y}^{\mydot})$ (and Grothendieck duality) to get an exact triangle,
$$
\xymatrix{ \omega_{E}^{\mydot} \ar[r] & \omega_{\overline X}^{\mydot}
  \ar[r] &
  \omega_{\tld X}^{\mydot} \tensor \OO_{\tld Y}(E) \ar[r]^-{+1} & \\
}
$$
Note that $\overline X$ is not equidimensional, but $\tld X$ and $E$ are (in fact,
they are Gorenstein and connected).  Since $E$ has dimension $n-1$ and $\tld X$ has
dimension $d$, we see that $h^{-n+1}(\omega_{\overline X}^{\mydot}) \cong
h^{-n+1}(\omega_E^{\mydot})$, proving the second statement.  Taking the $-d$th cohomology
proves the first statement since $\tld X$ and $E$ have normal crossings.  It is easy
to see that the third statement is true as well by taking any other cohomology.
These three facts will not be used directly, but they do suggest a way
to analyze $\myR \pi_* \omega_{\overline X}^{\mydot}$.

With this in mind, we now prove that we really only need to understand
the $-d$th cohomology of $\myR \pi_* \omega_{\overline X}^{\mydot}$,
at least in the Cohen-Macaulay case.

\begin{proposition}
  \label{CMSurjectiveImpliesDuBois}
  In addition to \eqref{setupNormal} assume further that $X$ is Cohen-Macaulay.  Then
  $X$ has DB singularities if and only if the natural map
  $$
  \myR^{-d} \pi_* \omega_{\overline X}^{\mydot} \rightarrow \omega_{X}
  $$
  is surjective (if and only if it is an isomorphism).
\end{proposition}

\begin{proof}
  If $X$ has DB singularities, the statement (including the one
  in parantheses) follows trivially from Lemma~\ref{Karldual}.

  Now assume that the natural map $\myR^{-d} \pi_* \omega_{\overline X}^{\mydot}
  \rightarrow \omega_{X}$ is surjective, but $X$ is not Du~Bois.
  Let $\Sigma_{DB}$ denote the non-Du~Bois locus of $X$ (cf.\
  \cite[2.1]{KovacsDuBoisLC1}) and let $x\in\Sigma_{DB}$ a general
  point of (a component of) $\Sigma_{DB}$. %
  By \cite[5.11]{SchwedeFInjectiveAreDuBois}, the natural map
  $$
  (\myR^i \pi_* \omega_{\overline X}^{\mydot})_x \rightarrow
  h^i(\omega_{X}^{\mydot})_x
  $$
  is injective for every $i$.  The right side of this equation is zero
  for $i \neq -d$ since $X$ is Cohen-Macaulay, and thus the left side
  is zero as well.  For $i = -d$, the map is surjective by assumption
  and, as we already noted, it is injective; hence it is an
  isomorphism.  In particular, the localized map $(\myR \pi_*
  \omega_{\overline X}^\mydot)_x \rightarrow (\omega^\mydot_{X})_x$ is
  a quasi-isomorphism,
  contradicting Lemma~\ref{Karldual} and the fact that $(X,x)$ is
  \emph{not} Du~Bois.
\end{proof}

\begin{remark}
  Alternatively, one could use general hyperplane sections to reduce to the case of
  an isolated non-Du~Bois point, and then apply local duality along with the key
  surjectivity of \cite{KovacsDuBoisLC2}.
\end{remark}

The following lemma will be important in the proof.

\begin{lemma}
\label{DimensionVanishingLemma}
Let $Z$ be a reduced closed subscheme of $Y$, a variety of finite type over $\bC$.
Then $h^i(\myR \sHom^{\mydot}_Y( \DuBois{Z}, \omega_Y^{\mydot})) = 0$ for $i < -\dim Z$.
\end{lemma}
\begin{proof}
  Without loss of generality, we may assume that $Z$ and $Y$ are affine.  Let $z \in
  Z$ be an arbitrary closed point. By local duality (see \cite[V, Theorem
  6.2]{HartshorneResidues} or \cite[2.4]{LipmanLocalCohomologyAndDuality}), it is
  sufficient to show that $\bH^j_z(Y, \DuBois{Z}) = \bH^j_z(Z, \DuBois{Z}) = 0$ for
  $j > \dim Z$. We consider the hypercohomology spectral sequence $H^p_z(Z,
  h^q(\DuBois{Z}))$ that computes this cohomology.  Note that $\dim (\Supp(
  h^q(\DuBois{Z}))) \leq \dim Z - q$ by \cite[V, 3.6]{GNPP}, so that $H^p_z(Z, h^q(\DuBois{Z})) = 0$ for $p
  > \dim Z - q$ (i.e., for $p + q > \dim Z$).  Therefore, we see that $\bH^j_z(Z,
  \DuBois{Z})$ vanishes for $j > \dim Z$ because every term in the spectral sequence
  that might possibly contribute to $\bH^j_z(Z, \DuBois{Z})$ is zero.
\end{proof}

\begin{corollary}\label{cor:vanishing-for-R}
  $\myR^i \pi_* \omega_E^{\mydot} = 0$ for $i < -\dim \Sigma$.
\end{corollary}

\begin{proof}
  By \cite{SchwedeEasyCharacterization} (cf.\ \eqref{setupNormal}) $\myR \pi_* \OO_E
  \qis \DuBois{\Sigma}$, so by Grothendieck duality
  $$
  \myR\pi_*\omega_E^\mydot\qis \myR\sHom_Y^{\mydot}(\DuBois{\Sigma},
  \omega_Y^\mydot).
  $$
  Then the statement follows from Lemma~\ref{DimensionVanishingLemma}.
\end{proof}

\begin{theorem}
\label{top}
If $X$ is equidimensional and $\codim_X\Sigma\geq 2$, then $\myR^{-d} \pi_* \omega_{\overline X}^{\mydot} \cong
\pi_* \omega_{\tld X}(E|_{\tld X})$.
\end{theorem}

\begin{remark}
  The codimension condition of this theorem implies that $X$ must be $R_1$. It is
  satisfied, for instance, if $X$ is normal and the maximal dimensional components of
  $\Sigma$ and $\Sing X$ coincide.
\end{remark}

\begin{proof}
  Applying $\myR \pi_*$ to the exact triangle,
  $$
  \omega_{E}^{\mydot} \rightarrow \omega_{\overline X}^{\mydot} \rightarrow
  \omega_{\tld X}^{\mydot} \tensor \OO_{\tld Y}(E) \rightarrow
  $$
  and taking cohomology leads to the exact sequence
  $$
  \myR^{-d} \pi_* \omega_{E}^{\mydot} \rightarrow \myR^{-d} \pi_* \omega_{\overline
    X}^{\mydot} \rightarrow \myR^{-d} \pi_* (\omega_{\tld X}^{\mydot} \tensor
  \OO_{\tld Y}(E)) \rightarrow \myR^{-d+1} \pi_* \omega_{E}^{\mydot}
  $$
  The outside terms are zero by Corollary~\ref{cor:vanishing-for-R} which implies that the middle two terms are isomorphic.  To complete the proof, simply observe that
\[
 \myR^{-d} \pi_* (\omega_{\tld X}^{\mydot} \tensor \OO_{\tld Y}(E)) = \myR^{-d} \pi_* (\omega_{\tld X}[d] \tensor \OO_{\tld Y}(E)) = \pi_* ( \omega_{\tld X}  \tensor \OO_{\tld Y}(E)) = \pi_* \omega_{\tld X}(E|_{\tld X}).
\]
\end{proof}

\begin{corollary}
  \label{CMNormalSurjectiveImpliesDuBois}
  In addition to \eqref{setupNormal} assume further that $X$ is normal and
  Cohen-Macaulay.  Then $X$ has DB singularities if and only if the natural map
  $$
  \pi_* \omega_{\tld X}(E|_{\tld X}) \to \omega_X
  $$
  (coming from~\eqref{CMSurjectiveImpliesDuBois} and~\eqref{top}) is surjective (if
  and only if it is an isomorphism).
\end{corollary}

\begin{remark}
  Note that $\pi_* \omega_{\tld X}(E|_{\tld X}) \to \omega_X$ is an isomorphism on
  the smooth locus of $X$ by construction, and hence it is always injective. For the
  same reason, this natural map is the same as the one coming from
  Lemma~\ref{PushdownInjects}.
\end{remark}

At this point, we have proven Theorem~\ref{Kempf} for a log resolution as in
\eqref{setupNormal}. The general case follows from the next Lemma which is well known
to experts (for example, it also follows from \cite[Lemma
1.6]{KawamataConeOfCurves}).


\begin{lemma}\label{lem:independent}
  Let $\pi : X' \rightarrow X$ be a proper birational morphism, $\Sigma\subseteq X$ a
  closed subset and denote by $G$ the reduced pre-image of $\Sigma$ via $\pi$. Assume
  that $\pi$ is chosen such that $X'$ is smooth and $G$ has simple normal crossings.
  Then $\pi_* \omega_{X'}(G)$ on $X$ is independent of the
  choice of $\pi$ up to natural isomorphism.
\end{lemma}

\begin{remark}
  This result is analogous to the fact that a multiplier ideal is independent of the
  resolution used to compute it.
\end{remark}

\begin{proof}
  Since any two proper birational morphisms mapping to $X$ with the required
  properties can be dominated by a third such, it is sufficient to prove the
  following: Let $X''$ be a smooth variety and $\phi: X''\to X'$ a proper birational
  morphism and let $H$ be the reduced pre-image of $\Sigma$ via $\pi\circ\phi$.
  Then $\phi_* \omega_{X''}(H) \simeq \omega_{X'}(G)$.

  The fact that $X'$ is smooth and $G$ is a simple normal crossing
  divisor implies that the pair $(X', G)$ has log canonical
  singularities.  Furthermore, the support of the strict transform of
  $G$ on $X''$ is contained in $H$ by definition, and the rest of the
  components of $H$ are $\phi$-exceptional. Therefore,
  $$
  \omega_{X''}(H)\simeq \phi^*(\omega_{X'}(G))(F)
  $$
  for an appropriate effective $\phi$-exceptional divisor $F$.
  Applying the projection formula yields
  $$
  \phi_*\omega_{X''}(H)\simeq \omega_{X'}(G)\otimes
  \phi_*\sO_{X''}(F),
  $$
  and since $F$ is effective and $\phi$-exceptional, $\phi_*\sO_{X''}(F)\simeq
  \sO_{X'}$ \cite[Lemma 1-3-2]{KMM}.
\end{proof}

We now turn our attention to using this critierion to prove that log canonical
singularities are Du~Bois.  First note that the statement is reasonably
straightforward if $X$ is Gorenstein so it would be tempting to try to take a
canonical cover, at least in the $\bQ$-Gorenstein case.  However, it is not clear
that the canonical cover of a Cohen-Macaulay log canonical singularity is also
Cohen-Macaulay.  Examples of rational singularities with non-Cohen Macaulay canonical
covers by \cite{SinghCyclicCoversOfRational} suggest that this might be too much to
hope for.  Therefore, a different technique will be used.

\newcommand{\tlpr}[1]{{#1}'}

\begin{lemma}
  \label{PushdownInjects}
  Let $X$ be a normal irreducible variety and let
  $\varrho: \tlpr X \rightarrow X$ be a log resolution of $X$.  Let $B$ be an
  effective integral divisor on $X$, $\tlpr B = \varrho^{-1}_* B$
  (the strict transform of $B$ on $\tlpr X$), and denote the reduced exceptional
  divisor of $\varrho$ by $G$.  Then there exists a natural injection,
  $$
  \varrho_* \omega_{\tlpr X}(\tlpr B+G) \into \omega_X(B).
  $$
\end{lemma}

\begin{proof}
  Let $\iota:U\into X$ denote the embedding of the open set over which $\varrho$ is
  an isomorphism. As $X$ is normal, we have that $\codim_X(X\setminus U)\geq 2$, and
  hence the following natural morphisms of sheaves:
  $$
  \varrho_* \omega_{\tlpr X}(\tlpr B+G)\into \left(\varrho_* \omega_{\tlpr
      X}(\tlpr B+G)\right)^{**}\simeq \iota_*\left(\varrho_* \omega_{\tlpr
      X}(\tlpr B+G)|_{U}\right) \simeq \omega_X(B),
  $$
  where $(\ )^{*}=\sHom_X(\__,\sO_X)$ denotes the dual of a sheaf. The isomorphisms
  follow because $X$ is $S_2$ and $\varrho$ is an isomorphism over $U$.
\end{proof}

\begin{lemma}
  \label{LCImpliesNicePushdown}
  Let $X$ be a normal irreducible variety with an effective $\bQ$-divisor $D$ such
  that $(X,D)$ has log canonical singularities, and let $\varrho: \tlpr X
  \rightarrow X$ be a log resolution of $(X, D)$.  Let $B$ be an effective integral
  divisor on $X$ with $B \leq D$.
  Denote the reduced exceptional divisor of $\varrho$ by $G$ and let $B' =
  \varrho^{-1}_* B$ and $\tlpr D = \varrho^{-1}_* D$.
  Then the following natural isomorphism holds:
  $$
  \varrho_* \omega_{\tlpr X}(B' + G) \simeq \omega_X(B)
  $$
\end{lemma}

\begin{proof}
  By Lemma~\ref{PushdownInjects} there exists a natural inclusion $\iota:\varrho_*
  \omega_{X'}(B'+G) \into \omega_X(B)$, so the question is local. We may assume that
  X is affine and need only prove that every section of $\omega_X(B)$ is already
  contained in $\varrho_*\omega_{X'}(B'+G)$. Note that $\iota$ restricted to the
  naturally embedded subsheaf $\omega_{X'}\subseteq \omega_{X'}(B'+G)$ gives the
  usual natural inclusion $\varrho_*\omega_{X'}\into \omega_X$.

  Next, choose a canonical divisor $K_{X'}$ and let $K_X=\varrho_*K_{X'}$. This is
  the divisor corresponding to the image of the section of $\omega_{X'}$
  corresponding to $K_{X'}$ via $\iota$.  As $D'=\varrho^{-1}_*D$, it follows that
  the divisors $K_{X'}+D'$ and $\varrho^{-1}_*(K_X+D)$ may only differ in exceptional
  components. We emphasize that these are actual divisors, not just equivalence
  classes (and so are $B$ and $B'$).

  Since $X$ and $X'$ are birationally equivalent, their function fields are
  isomorphic. Let us identify $K(X)$ and $K(X')$ via $\varrho^*$ and denote them by $K$.
  Further let $\sK$ and $\sK'$ denote the $K$-constant sheaves on $X$ and $X'$
  respectively.

  Now we have the following inclusions:
  $$
  \Gamma(X, \varrho_* \omega_{X'}(B'+G))\subseteq \Gamma(X,\omega_X(B)) \subseteq
  \Gamma(X,\sK)=K,
  $$
  and we need to prove that the first inclusion is actually an equality.  Let $g\in
  \Gamma(X,\omega_X(B))$.  By assumption, $B\leq D$, so
  \begin{equation}
    \label{eq:3}
      0\leq \divisor_X(g)+K_X+B\leq \divisor_X(g)+K_X +D
  \end{equation}
  As $(X,D)$ is log canonical and thus $K_X + D$ is $\bQ$-Cartier, there exists an $m\in\mathbb N$ such that $mK_X+mD$ is
  a Cartier divisor and hence can be pulled back to a Cartier divisor on $X'$.  By
  the choices we made earlier, we have that $\varrho^*(mK_X+mD)=mK_{X'}+mD'+\Theta$
  where $\Theta$ is an exceptional divisor. Again we emphasize that these are actual
  divisors and the two sides are actually equal, not just equivalent. We would also
  like to point out that $\Theta$ is not necessarily divisible by $m$, neither is
  $\varrho^*(mK_X+mD)$.

  However, using the fact that $(X,D)$ is log canonical, one obtains that $\Theta\leq
  mG$.  Combining this with (\ref{eq:3}) gives that
  $$
  0\leq \divisor_{X'}(g^m)+\varrho^*(mK_X+mD)\leq
  m\big(\divisor_{X'}(g)+K_{X'}+D'+G \big),
  $$
  and in particular we obtain that
  $$
  \divisor_{X'}(g) + K_{X'}+D'+G\geq 0.
  $$

  \begin{claim*}
    $\divisor_{X'}(g) + K_{X'}+B'+G\geq 0.$
  \end{claim*}

  \begin{proof}
    By construction
    \begin{equation}
      \label{eq:4}
      \divisor_{X'}(g) + K_{X'}+B'+G=\varrho^{-1}_*(\underbrace{\divisor_{X}(g) +
        K_{X}+B}_{\geq 0})+\underbrace{F+G}_{\text{exceptional}}.
    \end{equation}
    Where $F$ is an appropriate exceptional divisor, though it is not necessarily
    effective.  We also have that
    \begin{equation}
      \label{eq:5}
      \divisor_{X'}(g) + K_{X'}+B'+G=\underbrace{\divisor_{X'}(g) +
        K_{X'}+D'+G}_{\geq 0}-\underbrace{(D'-B')}_{\text{non-exceptional}}.
    \end{equation}
    Now let $A$ be an arbitrary irreducible component of $\divisor_{X'}(g) +
    K_{X'}+B'+G$. If $A$ were not effective, it would have to be exceptional by
    (\ref{eq:4}) and non-exceptional by (\ref{eq:5}). Hence $A$ must be effective and
    the claim is proven.
  \end{proof}
  \noindent
  It follows that
  $g\in\Gamma(X',\omega_{X'}(B'+G))=\Gamma(X,\varrho_*\omega_{X'}(B'+G))$, completing
  the proof.
\end{proof}

Now we are in a position to prove that Cohen Macaulay log canonical singularities are
Du~Bois.

\begin{theorem}
  \label{CMandLCImpliesDuBois}
  Suppose $X$ is normal and Cohen Macaulay and $\Delta \subset X$ an
  effective $\bQ$-divisor such that the pair $K_X + \Delta$ is
  $\bQ$-Cartier.  If $(X, \Delta)$ is log canonical, then $X$ has Du
  Bois singularities.
\end{theorem}

\begin{proof}
  Use Lemma~\ref{LCImpliesNicePushdown} setting $B = 0$ and note that the map of
  Corollary~\ref{CMNormalSurjectiveImpliesDuBois} is an isomorphism outside the
  singular locus.  Furthermore, both sheaves are reflexive (since they are abstractly
  isomorphic by Lemma~\ref{LCImpliesNicePushdown}), completing the proof.
\end{proof}

\section{The non-normal case}
\label{SectionNonnormalCase}

The aim of this section is to show that Cohen-Macaulay semi-log
canonical singularities are Du~Bois.  Let us begin by recalling some
of the relevant definitions.

\begin{definition}
  A reduced scheme $X$ of finite type over $\bC$ is said to be \emph{seminormal} if
  every finite morphism $X' \rightarrow X$ of reduced finite type schemes over $\bC$
  that is a bijection on points is an isomorphism.
\end{definition}

\begin{remark}
  If one is not working over an algebraically closed field of characteristic zero,
  one needs to alter the above definition somewhat.  See
  \cite{GrecoTraversoSeminormal} for details.
\end{remark}

\begin{definition}
  If $X$ is a reduced scheme of finite type over $\bC$ with
  normalization $\eta : X^N \rightarrow X$, then the \emph{conductor
    ideal sheaf of $X$ in its normalization} is defined to be the
  ideal sheaf $\Ann_{\OO_X}(\eta_* \OO_{X^N} / \OO_X)$.
\end{definition}

\begin{remark}
  \label{rem:affine-conductor}
  Consider the affine case where $R$ is a subring of its normalization $R^N$ (in its
  total field of fractions).  Then the conductor ideal is the largest ideal of $R^N$
  that is contained in $R$.  This implies that, with the previous notation, if $I_C$
  is the conductor ideal sheaf of $X$ in its normalization and if $I_B$ is the
  extension of $I_C$ to the normalization (that is $I_B = I_C \OO_{X^N}$), then
  $\eta_* I_B = I_C$.
\end{remark}

\begin{remark}\label{rem:codim-1-locus-of-sing}
  If $X$ is seminormal, then the conductor ideal sheaf of $X$ in its normalization is
  a radical ideal sheaf, even when extended to the normalization, see \cite[Lemma
  1.3]{TraversoPicardGroup}.  If $X$ is $S_2$, then all the associated primes of the
  conductor $I_C$ are height one, see \cite[Lemma 7.4]{GrecoTraversoSeminormal}, thus all
  the associated primes of $I_B$ are also height one (cf.\ \cite[9.2]{EisenbudCommutativeAlgebra}).
  Therefore, if $X$ is seminormal and $S_2$, then $\supp C$ is exactly the codimension
  $1$ locus of the singular set of $X$.
\end{remark}

\begin{definition}
  Let $X$ be a reduced equidimensional scheme of finite type over $\bC$. Assume that
  $X$ satisfies the following conditions:
  \begin{itemize}
  \item[(i)]  $X$ is $S_2$, and
  \item[(ii)] $X$ is seminormal.
  \end{itemize}
  These conditions imply that the conductor of $X$, in its normalization $X^N$, is a
  reduced ideal sheaf corresponding to an effective divisor on $X^N$ (cf.\
  \eqref{rem:codim-1-locus-of-sing}).  We let $B$ denote this divisor on $X^N$ and
  let $C$ denote the corresponding divisor on $X$ (by construction, these divisors have
  the same ideal sheaf in the normalization of $\OO_X$).
  Further, let $\Delta$ be an effective $\bQ$-divisor on $X$ and assume that $\Delta$
  and $C$ have no common components. Then $(X,\Delta)$ is said to be \emph{weakly
    semi-log canonical} if
  the pair $(X^N, B+\eta^{-1}_*\Delta)$ is log canonical.
\end{definition}

\begin{remark}\label{rem:conductor-reduced}
  The log canonical assumption implies that $B$ has to be reduced, but actually this
  does not impose any new conditions because simply from the fact that $X$ is
  seminormal, it follows by \cite[1.4]{GrecoTraversoSeminormal} that both $B$ and $C$
  are reduced.
\end{remark}

\begin{remark}
  The usual notion of \emph{semi-log canonical} also adds the
  additional condition
\begin{itemize}
\item[(iii-a)] $X$ is Gorenstein in codimension $1$,
\end{itemize}
This condition, under the previous seminormality assumption is equivalent to
\begin{itemize}
\item[(iii-b)] $X$ has simple double normal crossings in codimension
  1.
\end{itemize}
We won't need this condition, so we will leave it out.  Note that many
easy to construct schemes satisfy conditions (i) and (ii) but not
(iii).  For example, the reduced scheme consisting of the three axes
in $\bA^3$ does not have double crossings in codimension $1$, but is
both $S_2$ and seminormal.
\end{remark}

The key ingredient in the proof of the normal case is
Proposition~\ref{CMSurjectiveImpliesDuBois}, the injectivity of a certain map.  We
will now prove a strengthening of a special case of that injectivity that we will
need in this section.  First we need a lemma.

\begin{lemma}
  \label{ImprovedSupport}
  Let $d=\dim X$ and $\omegaplus{X}$ a complex that completes $\alpha$ to an exact
  triangle:
  $$
  \xymatrix{ \OO_X \ar[r]^{\alpha} & \DuBois{X} \ar[r] & \omegaplus{X} \ar[r]^-{+1}
    &}.
  $$
Then $\dim \Supp(h^i(\omegaplus{X})) \leq d - i - 1$ for
  $i\geq 0$ and $h^i(\omegaplus{X})=0$ for $i<0$.
\end{lemma}

\begin{proof}
  The statement follows from the definition for $i<0$, so we may assume that $i\geq
  0$.  Furthermore, for $i = 0$ the result also follows since $X$ is reduced and thus
  $X$ is generically smooth and hence generically Du~Bois, so $\OO_X \to
  h^0(\DuBois{X})$ is an isomorphism outside a set of codimension $1$.  We proceed by
  induction on the dimension of $X$.  Clearly it is true for zero (or even one)
  dimensional varieties (see \cite[4.9]{DuBoisMain} for the one-dimensional case).
  Let $\Sigma$ denote the singular set of $X$ and let $\pi : \tld X \rightarrow X$ be
  a resolution of $X$ coming from an embedded resolution as in \eqref{setupNormal}
  (so that $\pi$ is an isomorphism over the smooth locus of $X$).  Let $E =
  \pi^{-1}(\Sigma)_{\red}$, that is, $E$ is the reduced pre-image of the singular
  set.  We then have the exact triangle (cf.\ \eqref{setupNormal}),
  \begin{equation}
    \label{SupportTriangle}
    \xymatrix{
      \myR \pi_* \OO_{\tld X}(-E) \ar[r] & \DuBois{X} \ar[r] & \DuBois{\Sigma}
      \ar[r]^-{+1} &.
    }
  \end{equation}
  The case $i = 0$ follows from the fact that $h^0(\DuBois{X})$ is the structure
  sheaf of the seminormalization of $X$ by \cite[5.2]{SaitoMixedHodge} or
  \cite[5.6]{SchwedeFInjectiveAreDuBois}.  For $i > 0$, we note that it is sufficient
  to prove the statement for $h^i(\DuBois{X}) \cong h^i(\omegaplus{X})$. Observe that
  $$
  \dim \Supp(\myR^i \pi_* \OO_{\tld X}(-E)) \leq d - i -1\
  $$
  because $\pi|_{\tld X}$ is birational and the dimension of the exceptional set
  cannot be ``too'' big (cf.\ \cite[III.11.2]{Hartshorne}). By the inductive
  hypothesis the statement holds for $h^i(\DuBois{\Sigma})$ and so the long exact
  sequence coming from the triangle~\ref{SupportTriangle} completes the proof.
\end{proof}

Now we are in position to prove the desired injectivity statement.

\begin{proposition}
  \label{EasyInjection}
  In addition to \eqref{setupNormal}, let $\dim X=d$. Then the map $\myR^{-d} \pi_*
  \omega_{\overline X}^{\mydot} \rightarrow h^{-d}(\omega_X^{\mydot})$ is injective.
\end{proposition}

\begin{proof}
  By local duality, it is enough to show that $H^{d}_x(X, \OO_X) \rightarrow
  \bH^{d}_x(X, \DuBois{X})$ is surjective for every closed point $x \in X$.
  Considering the exact triangle,
  $$
  \xymatrix{ \OO_X \ar[r] & \DuBois{X} \ar[r] & \omegaplus{X} \ar[r]^-{+1} &}
  $$
  shows that it is enough to prove that $\bH^{d}_x(X, \omegaplus{X}) = 0$.  First
  observe that $H^p_x(X, h^q(\omegaplus{X})) = 0$ for $p > \dim
  \Supp(h^q(\omegaplus{X}))$.  Furthermore, by Lemma~\ref{ImprovedSupport}, we obtain
  that $H^p_x(X, h^q(\omegaplus{X})) = 0$ for $p > d - q - 1$, and therefore
  $\bH^d_x(X, \omegaplus{X}) = 0$ since for $p + q = d > d -1$ we see that every term
  in the standard spectral sequence that might contribute is already zero.
\end{proof}

The next proposition is the key step in our proof that Cohen-Macaulay semi-log
canonical singularities are Du~Bois.  This can be thought of as a generalization of
certain aspects of the Kempf-like criterion, Theorem~\ref{Kempf}.

We work in the following setting: $X$ is a reduced $d$-dimensional Cohen-Macaulay
scheme of finite type over $\bC$.  Let $\pi : \tld X \rightarrow X$ be a log
resolution of $X$, and $\Sigma\subseteq X$ a reduced closed subscheme such that $\pi$
is an iso\-mor\-phism outside $\Sigma$ (in particular $\Sigma\supseteq \Sing X$). Let
$F$ denote the reduced pre-image of $\Sigma$ in $\tld X$.  Consider the natural map
$\sI_{\Sigma} \rightarrow \myR \pi_* \OO_{\tld X}(-F)$ where $\sI_{\Sigma}$ is the ideal
sheaf of $\Sigma$.  Apply $\myR \sHom^{\mydot}_X(\blank, \omega_X^{\mydot} ) = \myR
\sHom^{\mydot}_X(\blank, \omega_X[d] ) $, which gives us a map
\[
\myR \sHom^{\mydot}_X(\myR \pi_* \OO_{\tld X}(-F), \omega_X[d] ) \to \myR
\sHom^{\mydot}_X(\sI_{\Sigma}, \omega_X[d] ).
\]
Then, by Grothendieck duality,
\[
\myR \sHom^{\mydot}_X(\myR \pi_* \OO_{\tld X}(-F), \omega_X[d] ) \cong \myR \pi_*
\myR \sHom^{\mydot}_{\tld X}(\OO_{\tld X}(-F), \omega_{\tld X}[d]) \cong \myR \pi_*
\omega_{\tld X}(F) [d].
\]
Taking the $-d$th cohomology gives us a natural map
\[
h^{-d}(\myR \pi_* \omega_{\tld X}(F) [d] ) \cong \pi_* \omega_{\tld X}(F) \to
\sHom_X(\sI_{\Sigma}, \omega_X ).
\]
We will use properties of this map to deduce that $X$ has DB singularities.

\begin{theorem}
  \label{keySlcDuBoisTheorem}
  Suppose we are in the setting described above.  If the natural map $\varrho : \pi_*
  \omega_{\tld X}(F) \rightarrow \sHom_X(\sI_{\Sigma}, \omega_X)$ is an isomorphism,
  then $X$ has DB singularities.
\end{theorem}

\begin{remark}
  \label{rem:F-is-B+G}
  Notice that if $X$ is not normal then $\Sigma$ contains the conductor.
\end{remark}

\begin{proof}
  By Lemma~\ref{lem:independent} the isomorphism class of $\pi_* \omega_{\tld X}(F)$
  is independent of the choice of the resolution, thus we may assume it came from an
  embedded resolution of $X$ in some $Y$ as in \eqref{setupNormal}.  In particular we
  have $\pi : \overline X \rightarrow X$, where $E$ is the reduced pre-image of
  $\Sigma$ in $\tld Y$ (i.e., $E|_{\tld X}=F$).  First, consider the following map of
  exact triangles,
  $$
  \xymatrix{%
    \sI_{\Sigma} \ar[r] \ar[d] & \OO_X \ar[d] \ar[r] & \OO_{\Sigma} \ar[d]
    \ar[r]^-{+1} & \\
    \myR \pi_* \OO_{\tld X}(-F 
    ) \ar[r] & \myR \pi_*
    \OO_{\overline X} \ar[r] & \myR \pi_* \OO_{E} \ar[r]^-{+1} & .
  }
  $$
  Applying $\myR \sHom_{X}^{\mydot}(\blank, \omega_X^{\mydot})$ produces
  $$
  \xymatrix{%
    \myR \pi_* \omega_E^{\mydot} \ar[r] \ar[d] & \myR \pi_* \omega_{\overline
      X}^{\mydot} \ar[r] \ar[d] & \myR \pi_* \omega_{\tld X}^{\mydot}(F 
    ) \ar[d] \ar[r]^-{+1} & \\
    \omega_{\Sigma}^{\mydot} \ar[r] & \omega_X^{\mydot} \ar[r] & \myR
    \sHom_{X}^{\mydot}(\sI_{\Sigma}, \omega_X^{\mydot}) \ar[r]^-{+1} &. \\
  }
  $$
 Considering the long exact cohomology sequence and using
  Corollary~\ref{cor:vanishing-for-R} leads to the following diagram:
  $$
  \xymatrix{%
    0 \ar[r] & h^{-d} (\myR\pi_* \omega_{\overline X}^{\mydot} ) \ar[r] \ar[d]_{\alpha} &
    h^{-d}( \myR \pi_* \omega_{\tld X}^{\mydot}(F 
    ) ) \ar[d]_{\beta} \ar[r] & h^{-d+1} (\myR \pi_*
    \omega_E^{\mydot} ) \ar[r] \ar[d]_{\gamma} & \ldots  \\
    0 \ar[r] & h^{-d}( \omega_X^{\mydot} \ar[r]) &
    h^{-d}(\myR \sHom^{\mydot}_X(\sI_{\Sigma}, \omega_X^{\mydot})) \ar[r] & h^{-d +
      1}(\omega_{\Sigma}^{\mydot}) \ar[r] & \\
  }
  $$
  Note that $\beta$ is simply the map $\varrho$ and thus it is surjective by
  hypothesis.  Note further that $\gamma$ is injective by Proposition~\ref{EasyInjection},
  thus $\alpha$ is surjective by the five lemma.  Combining this with
  Proposition~\ref{CMSurjectiveImpliesDuBois} completes the proof.
\end{proof}

We also need the following two lemmata.

\begin{notation}
  \label{not:notation-for-max-dim-component}
  Let $S$ be a reduced quasi-projective scheme of finite type of dimension $e$ over
  $\bC$. Then denote by $S_e$ the union of the $e$-dimensional irreducible components
  of $S$, and by $S_{<e}$ the union of the irreducible components of $S$ whose
  dimension is strictly less than $e$.
\end{notation}

\begin{lemma}
  \label{LemmaOnlyTopDimensionMatters}
  Let $\Sigma$ be a reduced quasi-projective scheme of finite type of dimension $e$
  over $\bC$.  Then $h^{-e}(\omega_{\Sigma}^{\mydot}) \simeq \omega_{\Sigma_e} =
  h^{-e}(\omega_{\Sigma_e}^{\mydot})$.
\end{lemma}

\begin{proof}
  Obviously, $\dim \Sigma_e = e$, $\dim \Sigma_{<e} < e$ and $\dim (\Sigma_e \cap
  \Sigma_{<e}) < e-1$.  Consider the short exact sequence
  $$
  0 \rightarrow \OO_{\Sigma} \rightarrow \OO_{\Sigma_e} \oplus \OO_{\Sigma_{<e}}
  \rightarrow \OO_{{\Sigma_e} \cap \Sigma_{<e}} \rightarrow 0
  $$
  where $\Sigma_e \cap \Sigma_{<e}$ is not necessarily a reduced scheme. Next, apply
  $\myR\sHom_\Sigma^{\mydot}(\__,\omega_\Sigma^\mydot)$ to get a long exact sequence:
  $$
  \ldots \rightarrow h^{-e}(\omega_{\Sigma_e \cap \Sigma_{<e}}^{\mydot}) \rightarrow
  h^{-e}(\omega_{\Sigma_e}^{\mydot} \oplus \omega_{\Sigma_{<e}}^{\mydot}) \rightarrow
  h^{-e}(\omega_{\Sigma}^{\mydot} ) \rightarrow h^{-e+1}(\omega_{\Sigma_e
    \cap \Sigma_{<e}}^{\mydot}) \rightarrow \ldots .
  $$
  As $\dim (\Sigma_e \cap \Sigma_{<e}) < e-1$, it follows that
  $h^{-e+1}(\omega_{\Sigma_e \cap \Sigma_{<e}}^{\mydot}) = h^{-e}(\omega_{\Sigma_e
    \cap \Sigma_{<e}}^{\mydot}) = 0$, and hence the statement holds.
\end{proof}


\begin{lemma}
  \label{LemmaOnIdentificationsOfHoms}
  Under the conditions of Theorem~\ref{keySlcDuBoisTheorem} and using the notation of
  \eqref{rem:codim-1-locus-of-sing} let $\eta : X^N \rightarrow X$ be the
  normalization of $X$ and assume that $\Sigma_e=C$ where $e=\dim\Sigma$.  Then
  $$
  \eta_* \omega_{X^N}(B) \simeq \eta_* \sHom_{X^N}(I_B, \omega_{X^N}) \simeq
  \sHom_{X}(I_C, \omega_X) \simeq \sHom_{X}(\sI_{\Sigma}, \omega_X).
  $$
\end{lemma}

\begin{proof}
   The first isomorphism holds because $\omega_{X^N}$ is a reflexive sheaf and $X^N$
  is normal.  The second follows from Grothendieck duality applied to the finite
  morphism $\eta$ (and the definition of the conductor).  To prove the last
  isomorphism, consider the map of short exact sequences
  $$
  \xymatrix{ 0 \ar[r] & \sI_{\Sigma} \ar[r] \ar[d] & \OO_{X} \ar[r]
    \ar@{=}[d] &
    \OO_\Sigma \ar[r] \ar[d] & 0\\
    0 \ar[r] & I_C \ar[r] & \OO_{X} \ar[r] & \OO_C \ar[r] & 0\\
  }
  $$
  and apply $\myR \sHom_X^{\mydot}(\blank, \omega_X^{\mydot})$. By
  Lemma~\ref{LemmaOnlyTopDimensionMatters}, we obtain
  $$
  h^{-d+1}(\myR \sHom^{\mydot}_X(\OO_\Sigma, \omega_X^{\mydot})) =
  h^{-d+1}(\omega_\Sigma^\mydot) \simeq
  h^{-d+1}(\omega_C^\mydot) =
  h^{-d+1}(\myR \sHom^{\mydot}_X(\OO_C, \omega_X^{\mydot})),
  $$
  thus we have the following diagram:
  $$
  \xymatrix{%
    0 \ar[r] \ar@{=}[d] & h^{-d}(\omega_X^{\mydot}) \ar[r] \ar@{=}[d] &
    h^{-d}(\myR \sHom^{\mydot}_X(I_{C}, \omega_X^{\mydot})) \ar[r] \ar[d] &
    h^{-d+1}(\omega_C^{\mydot}) \ar[r] \ar[d]^\simeq & h^{-d+1}(\omega_X^{\mydot})
    \ar@{=}[d] \\
    0 \ar[r] & h^{-d}(\omega_X^{\mydot}) \ar[r] & h^{-d}(\myR \sHom^{\mydot}_X(\sI_{\Sigma},
    \omega_X^{\mydot})) \ar[r] & h^{-d+1}(\omega_{\Sigma}^{\mydot}) \ar[r] &
    h^{-d+1}(\omega_X^{\mydot}).
  }
  $$
  The statement then follows from the five lemma.
\end{proof}

\begin{theorem}
  \label{CMslcImpliesDuBois}
  If $(X,\Delta)$ is a Cohen-Macaulay weakly semi-log canonical pair, then $X$ has
  DB singularities.
\end{theorem}

\begin{proof}
  In our setting we may assume that $X$ is non-normal but that it is $S_2$ and
  seminormal.  Thus $C$, the non-normal locus of $X$, is a codimension $1$ subset of
  $X$ (see Remark \ref{rem:codim-1-locus-of-sing}), in particular $C \neq 0$.  Using
  the same notation as above, let $\Sigma=\Sing X$.  Also observe that any log
  resolution $\pi:\tld X\to X$ factors through $\eta$:
  $$
  \xymatrix{ \tld X \ar[r]^{\zeta} \ar@/_1pc/[rr]_\pi & X^N \ar[r]^{\eta} & X.  }
  $$
  Therefore $\pi_* \omega_{\tld X}(F)= \eta_*\zeta_*\omega_{\tld X}(F)$. Recall from
  Remark~\ref{rem:conductor-reduced} that $C$ is reduced, so $\Sigma_e=C$, and then
  by Lemma~\ref{LemmaOnIdentificationsOfHoms}, $\eta_* \omega_{X^N}(B) \simeq
  \sHom_{X}(\sI_{\Sigma}, \omega_X)$. Then by Theorem~\ref{keySlcDuBoisTheorem}
  it is sufficient to show that $\zeta_* \omega_{\tld X}(F) = \omega_{X^N}(B)$.
  Notice that by definition $F=\zeta^{-1}_*B+G$ where $G$ is the reduced exceptional
  divisor of $\zeta$ (cf.\ Remark~\ref{rem:F-is-B+G}). By assumption $(X^N,
  B+\eta^{-1}_*\Delta)$ is log canonical, so Lemma~\ref{LCImpliesNicePushdown}
  implies that $\zeta_* \omega_{\tld X}(F) = \omega_{X^N}(B)$.
\end{proof}

\section{Cohomologically insignificant degenerations}

DB singularities were originally defined by Steenbrink in a Hodge-theoretic
context and they admit many interesting properties. In particular, they are strongly
connected with cohomologically insignificant degenerations.

Inspired by Mumford's definition of insignificant surface singularities Dolgachev
\cite{DolgachevCohomologicallyInsignificant} defined a \emph{cohomologically
  insignificant degeneration} as follows:

\newcommand\XX{\scr X}%
\newcommand{\spit}{{\rm sp}^j_t}%
Let $f : \XX\to S$ be a proper holomorphic map from a complex space $\XX$ to the unit
disk $S$ that is smooth over the punctured disk $S \setminus \{0\}$.  For $t\in S$
denote $\XX_t = f^{-1}(t)$. The fiber $\XX_0$, the \emph{special fiber}, may be
considered as a degeneration of any fiber $\XX_t$, $t\neq 0$.  Let $\beta_t :
H^j(\XX)\to H^j(\XX_t)$ be the restriction map of the $j^\text{th}$-cohomology spaces
with real coefficients. Because $\XX_0$ is a strong deformation retract of $\XX$ the
map $\beta_0$ is bijective. The composite map
$$
\spit=\beta_t\circ\beta_0^{-1}: H^j(\XX_0)\to H^j(\XX_t)
$$
is called the \emph{specialization map} and plays an important role in the theory of
degenerations of algebraic varieties. According to Deligne for every complex
algebraic variety $Y$ the cohomology space $H^n(Y)$ admits a canonical and functorial
mixed Hodge structure.  However in general $\spit$ is not a morphism of these mixed
Hodge structures. On the other hand, Schmid \cite{SchmidVariationOfHodge} and Steenbrink
\cite{SteenbrinkLimitsOfHodge} introduced a mixed Hodge structure on $H^i(\XX_0)$, the
\emph{limit Hodge structure}, with respect to which $\spit$ becomes a morphism of
mixed Hodge structures.

In the above setup, $\XX_0$ is called a \emph{cohomologically $j$-insignificant
degeneration} if $\spit$ induces an isomorphism of the $(p, q)$-components of the
mixed Hodge structures with $pq = 0$. Note that this definition is independent of the
choice of $t\neq 0$. Finally, $\XX_0$ is called a \emph{cohomologically
  insignificant degeneration} if it is cohomologically $j$-insignificant for every
$j$.

For us the relevance of this notion is that Steenbrink
\cite{SteenbrinkCohomologicallyInsignificant} proved that every proper, flat
degeneration $f$ over the unit disk $S$ is cohomologically insignificant provided
$f^{- 1}(0)$ has DB singularities.  As a combination of Steenbrink's result and
our main theorem, we obtain the following.

\begin{theorem}
  \label{CMslcImpliesInsignificant}
  Let $X$ be a proper algebraic variety over $\bC$ with Cohen-Macaulay semi-log
  canonical singularities.  Then every proper flat degeneration $f$ over the unit
  disk with $f^{-1}(0) = X$ is cohomologically insignificant.
\end{theorem}

\section{Kodaira vanishing}

Following ideas of Koll\'ar \cite[\S 9]{KollarShafarevich}, we give a proof of
Kodaira vanishing for log canonical singularities.  This was recently also proven by
Fujino using a different technique \cite[Corollary
5.11]{FujinoVanishingAndInjectivity}.  Our proof applies to (weakly) semi-log canonical
singularities as well, while Fujino's does not use the Cohen-Macaulay assumption.
Partial results were also obtained by Koll\'ar \cite[12.10]{KollarShafarevich} and
Kov\'acs \cite[2.2]{KovacsDuBoisLC2}.

\begin{convention}
  For the rest of the section, all cohomologies are in the Euclidean topology.
  Nonetheless, the results remain true for coherent cohomology by Serre's GAGA
  principle.
\end{convention}

First we need a variation on an important theorem.

\begin{theorem}[{\cite[Theorem 9.12]{KollarShafarevich}}]
  \label{thm:several}
  Let $X$ be a proper variety and $\sL$ a line bundle on $X$. Let $\sL^n\simeq \OO_X
  (D)$, where $D=\sum d_i D_i$ is an effective divisor (the $D_i$ are the irreducible
  components of $D$).  We also assume that the generic point of each $D_i$ is a smooth point of $X$.  Let $s$ be a global section of $\sL^n$ whose zero divisor is $D$.
  Assume that $0<d_i<n$ for every $i$. Let $Z$ be the scheme obtained by taking the
  $n^{\text{th}}$ root of $s$ (that is, $Z = X[\sqrt{s}]$ using the notation from
  \cite[9.4]{KollarShafarevich}). Assume further that
  \begin{equation*}
    H^j(Z, \bC_Z) \to H^j(Z, \OO_Z)
  \end{equation*}
  is surjective. Then for any collection of $b_i\geq 0$ the natural map
  $$
  H^j\left(X, \sL^{-1}\left(-\sum b_iD_i\right)\right) \to H^j(X, \sL^{-1})
  $$
  is surjective.
\end{theorem}

The aforementioned variation of this theorem differs from the version stated above in that
$Z$ was defined to be the normalization of $X[\sqrt{s}]$. However, as we are dealing with
possibly non-normal schemes (e.g., \emph{slc}) we need this version. In some
sense, the proof of this formulation is actually easier and we sketch the argument
below.  The strategy is the same as in \cite[Theorem 9.12]{KollarShafarevich}, but as
$Z$ is not normalized, some steps and ingredients are different.

\begin{proof}
  Let $Z=X[\sqrt s]=\operatorname{Spec}_X \sum_{t=0}^{n-1}\sL^{-t}$ and $p : Z
  \rightarrow X$ denote the natural map. By construction there is a decomposition
  $p_* \sO_Z \simeq \sum_{t = 0}^{n-1} \sL^{-t}$.  Fixing an $n^{\text{th}}$ root of unity
  $\zeta$ and considering the associated $\bZ/n$-action on $Z$ (coming from the
  cyclic cover), we see that $\bZ/n$ acts on the summand $\sL^{-t}$ by multiplication
  by $\zeta^{-t}$, so $p_* \sO_Z \simeq \sum \sL^{-t}$ is actually the eigensheaf
  decomposition of $p_*\sO_Z$ cf.\ \cite[Proposition 9.8]{KollarShafarevich}.  One
  also has an eigensheaf decomposition, $p_* \bC_Z \simeq \sum_{t = 0}^{n-1} \sG_t$.
  For ease of reference, set $\sG_t$ to be the eigensheaf corresponding to the
  eigenvalue $\zeta^{-t}$. In particular, $\sG_t\subseteq \sL^{-t}$.


  With these decompositions, note that we have a surjective map
  $$
  \sum_t H^j(X, \sG_t) \simeq H^j(Z,\bC_Z) \twoheadrightarrow H^j(Z,\sO_Z)\simeq
  \sum_t H^j(X, \sL^{-t})
  $$
  and so in particular we have a surjection
  \begin{equation}
    \label{eq:2}
    H^j(X, \sG_1) \twoheadrightarrow H^j(X, \sL^{-1})
  \end{equation}
  for every $j$.  We now claim that $\sG_1$ is a subsheaf of $\sL^{-1}(-\Sigma b_i
  D_i)$. As they are both subsheaves of $\sL^{-1}$, this is a local question and it
  is enough to show that for every connected open set $U\subseteq X$, the inclusion
  $\gamma:\Gamma(U,\sG_1)\into \Gamma(U, \sL^{-1})$ factors through the inclusion
  $\delta:\Gamma(U,\sL^{-1}((-\Sigma b_i D_i))\into \Gamma(U, \sL^{-1})$.

  If $U$ is such that $U\cap D=\emptyset$, then $\delta$ is an isomorphism and so the
  statement holds trivially.

  If $U$ is such that $U\cap D\neq\emptyset$, then $\Gamma(U, \sG_1) = 0$ by
  \eqref{claim:no-sections-of-G1}, and so the statement again holds trivially.

  \begin{claim}\label{claim:no-sections-of-G1}
    Let $U\subseteq X$ be a connected open set such that $U\cap D\neq\emptyset$. Then
    $\Gamma(U, \sG_1) = 0$.
  \end{claim}

  \begin{proof}
    We give two short proofs of this claim.

    First, let $\tld Z$ denote the normalization of $Z$, $\tld p:\tld Z\to X$ the
    induced map, and $\tld \sG_t$ the eigensheaf of the $\bZ/n$ action on $\tld
    p_*\bC_{\tld Z}$ corresponding to the eigenvalue $\zeta^{-t}$. Then
    $p_*\bC_Z\into \tld p_*\bC_{\tld Z}$ naturally, so in particular, $\sG_t\subseteq
    \tld\sG_t$ for all $t$ and hence the statement follows by
    \cite[9.11.3]{KollarShafarevich}.

    Alternatively, one can give a direct proof as follows. The assumptions imply that
    there exists a dense open subset $U'\subseteq U$ such that each $x\in U'\cap D$
    has a neighborhood where $X$ is smooth and $D$ is defined by a power of a
    coordinate function. Then the computation in \cite[9.9]{KollarShafarevich} shows
    that $p^{-1}U'$ and therefore $p^{-1}U$ are connected and the claim follows
    easily.
  \end{proof}
  Therefore $\sG_1$ is indeed a subsheaf of $\sL^{-1}(-\Sigma b_i D_i)$ and one
  obtains a factorization
  $$
  \xymatrix{%
    H^j(X, \sG_1) \ar@/^1pc/@{>>}[rr] \ar[r] & H^j\left(X, \sL^{-1}\left(-\sum
        b_iD_i\right)\right) \ar[r] & H^j(X, \sL^{-1}), }
  $$
  where the composition is surjective by (\ref{eq:2}) and so the second arrow is
  surjective as well.
\end{proof}

\begin{theorem}[(Serre's vanishing)] \cite[5.72]{KollarMori}
  \label{thm:cm}
  Let $X$ be a projective scheme over a field $k$ of pure dimension $n$ with ample
  Cartier divisor $D$. Then the following are equivalent:
  \begin{enumerate}
  \item $X$ is Cohen-Macaulay,
  \item $H^j(X,\OO_X(-rD))=0$ for every $j<n$ and $r\gg 0$.
  \end{enumerate}
\end{theorem}

In order to use the ``usual" covering trick, we need to establish that our
assumptions are inherited by the covers. Examples of rational singularities with
non-Cohen Macaulay canonical covers in \cite{SinghCyclicCoversOfRational} suggest
that this is not entirely obvious.

First we need the following construction.

\begin{notation}\rm
  Let $\tau:S\to T$ be a finite morphism between reduced schemes of finite type over
  $\bC$ and assume that $T$ and $S$ are normal. Let $i:U=T\setminus\Sing T\into T$ and
  $j:V=\tau^{-1}U\into S$.  Since $T$ is normal and $\tau$ is finite, $\codim_T (T \setminus U)\geq
  2$ and $\codim_S (S \setminus V)\geq 2$.  Let $D\subseteq T$ be an effective Weil divisor.  Then
  $D|_U$ is a Cartier divisor corresponding to the invertible sheaf $\sL$ on $U$ and
  $D$ induces a section of $\sL$: $\delta:\sO_U\into \sL$. Furthermore,
  $\sO_T(D)\simeq i_*\sL$.  Then the pullback of $\delta$ induces a section
  $\tau^*\delta:\sO_V\into \tau^*\sL$, which in turn induces a section of the rank
  $1$ reflexive sheaf $j_*\tau^*\sL$ on $S$. Denote the corresponding Weil divisor by
  $\tau^{[*]}D$, i.e., $\sO_S(\tau^{[*]}D)\simeq j_*\tau^*\sL$.
  An alternative way to obtain $\tau^{[*]}D$ is to take the closure (in $S$) of the
  divisor $\tau^*(D|_U)\subseteq V$.
\end{notation}

\begin{lemma}\label{lem:cover-is-cm}
  Let $(X,\Delta)$ be a projective log variety with weakly semi-log canonical
  singularities and $\sigma: Z\to X$ a cyclic cover of $X$ induced by a general
  section of a sufficiently large power of an ample line bundle as in \cite[2.50]{KollarMori}.  Then
  there exists an effective $\bQ$-divisor $\Gamma$ on $Z$ such that $(Z,\Gamma)$ is
  weakly semi-log canonical.
  Furthermore, if $X$ is Cohen-Macaulay, then so is $Z$. 
\end{lemma}

\begin{proof}
  Let $\sL$ be an ample line bundle on $X$, $m\gg 0$ and $s\in \coh 0.X.\sL^m.$ a
  general section. Then $D=(s=0)$ is reduced.  As before, let $\eta:X^N\to X$ denote
  the normalization of $X$ and $B\subset X^N$ the extension of the conductor to $X^N$
  (cf.\ \eqref{rem:codim-1-locus-of-sing}). Then, by assumption, $(X^N,
  B+\eta^{-1}_*\Delta)$ is log canonical. Observe that as $\eta$ is finite,
  $\eta^*\sL$ is also ample. Therefore, by \cite[5.17.(2)]{KollarMori}, $(X^N,
  B+\eta^{-1}_*\Delta+\eta^*D)$ is also log canonical.

  Let $\sA=\oplus_{i=0}^{m-1}\sL^{-i}$ with the $\ring X.$-algebra structure induced
  by $s$. Let $\sigma:Z=\operatorname{Spec}_X \sA \to X$ be as in
  \cite[2.50]{KollarMori} and similarly $\tld\sigma:W=\operatorname{Spec}_X \eta^*\sA
  \to X^N$.

  $$
  \xymatrix{%
    W \ar[r]^{\tld\sigma} \ar[d]_\tau & {X^N\hskip-8pt} \ar[d]^\eta \\
    Z \ar[r]_\sigma & X
  }
  $$

  By assumption $X^N$ is normal, i.e., $R_1$ and $S_2$ by Serre's criterion. Then $W$
  is also $R_1$ by \cite[2.51]{KollarMori} and furthermore $\tld\sigma_*\ring
  W.\simeq \eta^*\sA$ is also $S_2$ (as it is locally free), so we see that $W$ is
  also $S_2$ by \cite[5.4]{KollarMori}. Therefore $W$ is normal by Serre's criterion
  and hence $\tau$ factors through the normalization of $Z$:
  $$
  \xymatrix{%
    W \ar[r]_{\widetilde\tau} \ar@/^1pc/[rr]^\tau & {\widetilde Z} \ar[r] & Z
  }
  $$
  From the construction it is clear that $\tau:W\to Z$ is birational, and hence so is
  $\widetilde\tau$. However, then it must be an isomorphism by Zariski's Main
  Theorem, and hence $W$ is the normalization of $Z$.

  As $X$ is $S_2$, so is $\sigma_*\sO_Z\simeq \sA$, and then $Z$ is $S_2$ as well.

  Let $T=(\Sing X\cap D)_\text{red}\subseteq X$ and $\tld T=\sigma^{-1}T\subseteq Z$
  closed subsets in $X$ and $Z$ respectively. Then by construction $\codim_XT\geq 2$
  and hence $\codim_Z\tld T\geq 2$. Observe that for any $z\in Z\setminus \tld T$
  either $Z$ is smooth at $z$ or $\sigma$ is \'etale in a neighborhood of $z$ in $Z$.
  Therefore, $Z$ is seminormal in codimension $1$ (since $X$ is seminormal) and as we
  have just shown that $Z$ is also $S_2$, it follows by
  \cite[2.7]{GrecoTraversoSeminormal} that $Z$ is seminormal everywhere.

  Let $M=(\tld\sigma^*(\eta^* D))_{\text{red}}$.  Since $X^N$ and $W$ are smooth in
  codimension $1$, and $\tld\sigma$ is ramified exactly along $\eta^*D$ with degree
  $m$ and ramification index $m$ everywhere, it follows, that
  $$
  \tld\sigma^*(K_{X^N}+ B+\eta^{-1}_*\Delta+\eta^*D)= K_W+\tld\sigma^{[*]}(
  B+\eta^{-1}_*\Delta)+ M.
  $$
  Then $(W,\tld\sigma^{[*]}( B+\eta^{-1}_*\Delta)+M)$ is log canonical by
  \cite[5.20(4)]{KollarMori}.

  Let $B_Z$ denote the subscheme defined by the conductor of $Z$ in $W$.  We claim
  that $B_Z$ is contained in the proper transform of $B$, i.e., $B_Z\subseteq
  \tld\sigma^{[*]}B$ (cf.\ \eqref{rem:affine-conductor}). To see this, note that
  since $Z$ is $S_2$ and seminormal, the conductor is simply the codimension $1$ part
  of the non-smooth locus of $Z$.  But the non-smooth locus of $Z$ is the pre-image
  of the non-smooth locus of $X$ by \cite[2.51]{KollarMori} and so the claim follows.

  Then there exists an effective $\bQ$-divisor $\Theta$ on $W$ satisfying
  $\tld\sigma^{[*]}(B+\eta^{-1}_*\Delta)=B_Z+\Theta$.  Now, if we choose
  $\Gamma=\tau_*(\Theta+M)$, then $(Z,\Gamma)$ has weakly semi-log canonical
  singularities.

  If $X$ is Cohen-Macaulay, then so is $\sA\simeq \pi_*\OO_Z$ and $Z$ is
  Cohen-Macaulay by \cite[5.4]{KollarMori}.
\end{proof}

\begin{corollary}
  \label{KodairaVanishingForLC}
  Kodaira vanishing holds for Cohen-Macaulay weakly semi-log canonical varieties: Let
  $(X,\Delta)$ be a projective Cohen-Macaulay weakly semi-log canonical pair and
  $\sL$ an ample line bundle on $X$. Then $\coh i.X.\sL^{-1}.=0$ for $i<\dim X$.
\end{corollary}

\begin{proof}
  We will use the notation from \eqref{lem:cover-is-cm}. $Z$ is Cohen-Macaulay and
  $(Z,\Gamma)$ is weakly semi-log canonical. Therefore by \eqref{CMslcImpliesDuBois}
  $Z$ is Du~Bois, and it follows from \eqref{thm:several} that $\coh i.X.\sL^{-m}.\to
  \coh i.X.\sL^{-1}.$ is surjective for all $m\geq 0$. Serre's vanishing
  \eqref{thm:cm} implies that $\coh i.X.\sL^{-m}.=0$ for $m\gg 0$ and $i<\dim X$, so
  the desired statement follows.
\end{proof}

This implies invariance of plurigenera in stable Gorenstein families.

\begin{corollary}
  Let $f:\XX\to S$ be a stable Gorenstein family, i.e., a flat projective family of
  canonically polarized varieties with at most Gorenstein (weakly) semi-log canonical
  singularities. Then $h^i(\XX_t,\omega_{\XX_t}^{m})$ is independent of $t\in S$ for
  any $m>0$ and $i\geq 0$.
\end{corollary}

\begin{proof}
  By Serre duality
  $h^i(\XX_t,\omega_{\XX_t}^{m})=h^{\dim\XX_t-i}(\XX_t,\omega_{\XX_t}^{-m+1})$ and
  the latter vanishes for $i>0$ and $m>1$ by \eqref{KodairaVanishingForLC}. Then
  $$
  h^0(\XX_t,\omega_{\XX_t}^{m})=\chi(\XX_t,\omega_{\XX_t}^{m})
  $$
  for $m>1$ and this is independent of $t$ because $f$ is flat.

  So we may assume that $m=1$. Then
  $h^i(\XX_t,\omega_{\XX_t})=h^{\dim\XX_t-i}(\XX_t,\sO_{\XX_t})$ and this is
  independent of $t$ for any $i$ by \eqref{CMslcImpliesDuBois} and \cite[Th\'eor\`eme
  4.6]{DuBoisMain}.
\end{proof}



\providecommand{\bysame}{\leavevmode\hbox to3em{\hrulefill}\thinspace}
\providecommand{\MR}{\relax\ifhmode\unskip\space\fi MR}
\providecommand{\MRhref}[2]{%
  \href{http://www.ams.org/mathscinet-getitem?mr=#1}{#2}
}
\providecommand{\href}[2]{#2}

\end{document}